\documentclass{article}
\usepackage{indentfirst}
\usepackage{amsmath}
\usepackage{amsthm}
\usepackage{amsfonts}
\usepackage{graphics}
\usepackage{mathrsfs}
\usepackage{latexsym}
\usepackage{flafter}
\usepackage{amssymb}
\usepackage{dsfont}
\usepackage{colortbl}
\usepackage{cite}
\allowdisplaybreaks[4]
\newtheorem{theorem}{Theorem}

\newtheorem{lemma}{Lemma}

\begin{document}
\title{\textbf{Explicit error bound of the fast multipole method for scattering problems in 2-D}}
\author{Wenhui Meng\thanks{E-mail address: mwh@nwu.edu.cn. This paper has been accepted by Calcolo.}\\
\small\parbox{120mm}{\begin{center}\emph{School of Mathematics, Northwest University, Xi'an, $710127$, China}
\end{center}}}
\date{}
\maketitle

\begin{abstract}
This paper is concerned with the error estimation of the fast multipole method (FMM) for scattering problems in 2-D. The FMM error is caused by truncating Graf's addition theorem in each step of the algorithm, including two expansions and three translations. We first give a novel bound on the truncation error of Graf's addition theorem by the limiting forms of Bessel and Neumann functions, and then estimate the error of the FMM. Explicit error bound and its convergence order are derived. The method proposed in this paper can also be used to the FMM for other problems, such as potential problems, elastostatic problems, Stokes flow problems and so on.

\end{abstract}

\section{Introduction}
The fast multipole method (FMM) proposed by Rokhlin \cite{1} that has been widely applied in solving particle interaction problems and boundary integral equations. For solving a dense linear system with $N$ unknowns by an iterative method, it will require $\mathcal{O}(N^2)$ operations to store the matrix and compute the matrix-vector product. FMM can reduce the computing time and memory requirement to $\mathcal{O}(N)$ or $\mathcal{O}(N\ln{N})$. In recent years, FMM has been widely used to solve many mathematical physics problems, such as potential problems, elastostatic problems, acoustic and electromagnetic scattering problems. Some applications of the FMM for solving the scattering problems can be found in Refs. \cite{2,3,4,5,6,7,8,9}.

The error of FMM is caused by truncating the expansions and translations in the algorithm, including the multipole and local expansions, M2M, M2L and L2L translations. Most of the existing work focused on estimating the truncation errors of the multipole and local expansions. In Refs. \cite{2,3,4,5}, the truncation errors of the multipole and local expansions about the FMM for scattering problems in 2-D were studied. Some estimates about the expansion errors of the FMM in 3-D were proposed in Refs. \cite{9,10,11}. In Refs. \cite{6,9,14}, some empirical formulas were proposed to determine the truncation number of the multipole and local expansions.

There is little existing work studying the truncation errors of the M2M, M2L and L2L translations, possibly because the formulas used in these translations are more complex than that used in the multipole and local expansions. The article by Christophe \cite{12} estimated the error of the M2L translation of the FMM for scattering problems. However, the estimation has a very stringent assumption, which leads to its weak applicability. The article by Zhang, Wang and Cai \cite{13} proved the exponential convergence of the expansions and translations in the FMM for scattering problems in layered media. In fact, the error of FMM is the superposition of expansion errors and translation errors. In the FMM for solving 2-D potential problems, Wala and Kl$\ddot{\mathrm{o}}$ckner \cite{20} estimated the error of a translation chain: source$\rightarrow$multipole$\rightarrow$local$\rightarrow$local, and derived an explicit bound on the error.

The aim of this paper is to analyze the error of the FMM for 2-D scattering problems, including expansion errors, translation errors and their superposition, i.e., the error of a translation chain: source$\rightarrow$multipole$\rightarrow$multipole$\rightarrow$local$\rightarrow$local. The FMM formulations for 2-D scattering problems are based on Graf's addition theorems \cite{17,18} of Bessel function $J_n(z)$ and Hankel function $H^{(1)}_n(z)$. Thus, the key issue lies in how to give sharp estimates of the truncation errors of the addition theorems. Amini and Profit \cite{5} gave an estimate of the truncation error, but the result contains a higher order derivative. A more explicit bound on the truncation error was proposed in \cite{19}. However, those results are based on the asymptotic behavior of $J_n(z)$ and $H^{(1)}_n(z)$ as $n\rightarrow\infty$. This leads to the bounds which are not sharp for smaller truncation numbers, and then can not be well used to estimate the error of FMM.

In this paper, the bounds on truncation errors of Graf's addition theorems will be estimated by using the limiting forms of Bessel functions (when $z\rightarrow0$), and the results will be further applied to the error estimation of FMM. Section 2 gives a quick general overview of the FMM error for scattering problems. In Section 3, we propose novel bounds on truncation errors of Graf's addition theorems. Section 4 gives an explicit bound on the FMM error. We apply the previous result to the adaptive tree structure, convergence order of the FMM is obtained.

\section{Error of the FMM for scattering problems}
Consider the time-harmonic acoustic or electromagnetic waves propagating in homogeneous medium, and impinging on an infinitely long cylindrical obstacle. Denote the cross section of the obstacle by $\Omega$ and assume that $\Omega\in\mathds{R}^2$ is a bounded domain with a connected boundary $\partial\Omega$. Thus, the scattered wave $u(\mathbf{x})$ satisfies the Helmholtz equation
$$\triangle u(\mathbf{x})+k^2u(\mathbf{x})=0,\qquad\mathbf{x}\in{\mathds{R}^2\backslash{\overline{\Omega}}},$$
and the Sommerfeld radiation condition
$$\lim_{|\mathbf{x}|\rightarrow\infty}\sqrt{|\mathbf{x}|}\left(\frac{\partial{u}(\mathbf{x})}{\partial{|\mathbf{x}|}}-\mathrm{i}ku(\mathbf{x})\right)=0,$$
where $k>0$ is the wave number. The solution of the above problem can be represented as the form of single- or double-layer potential \cite{15}:
$$(S\varphi)(\mathbf{x})=\int_{\partial\Omega}\Phi(\mathbf{x},\mathbf{y})\varphi(\mathbf{y})\mathrm{d}s(\mathbf{y}),\qquad\mathbf{x}\in{\mathds{R}^2\backslash{\partial\Omega}},$$
$$(K\varphi)(\mathbf{x})=\int_{\partial\Omega}\frac{\partial\Phi(\mathbf{x},\mathbf{y})}{\partial\nu(\mathbf{y})}\varphi(\mathbf{y})\mathrm{d}s(\mathbf{y}),
\qquad\mathbf{x}\in{\mathds{R}^2\backslash{\partial\Omega}},$$
where the density $\varphi$ is an integrable function, $\Phi(\mathbf{x},\mathbf{y})$ is the fundamental solution to the Helmholtz equation, in 2-D which is given by
$$\Phi(\mathbf{x},\mathbf{y})=\frac{\mathrm{i}}{4}H^{(1)}_0(k|\mathbf{x}-\mathbf{y}|),\qquad \mathbf{x}\neq\mathbf{y}.$$
where i is the imaginary unit. Let the solution satisfy the corresponding boundary conditions, the boundary integral equations about $(S\varphi)(\mathbf{x})$ and $(K\varphi)(\mathbf{x})$ are obtained.

Assume that $\partial\Omega$ is divided into elements $\Delta\partial\Omega_j(j=1,\ldots,N)$, and place one node on each element. Thus, the single-layer potential has the following discretized form:
$$(S\varphi)(\mathbf{x})=\frac{\mathrm{i}}{4}\sum_{j=1}^{N}\int_{\Delta\partial\Omega_j}H^{(1)}_0(k|\mathbf{x}-\mathbf{y}_j|)\varphi(\mathbf{y}_j)\mathrm{d}s(\mathbf{y}_j),\qquad\mathbf{x}\in\partial\Omega,$$
where $\mathbf{y}_j$ is the node on $\Delta\partial\Omega_j$. For all nodes $\mathbf{x}_i\in\Delta\partial\Omega_i(i=1,\ldots,N)$, the above discretized form can be written as $\mathbf{A}\mathbf{z}$, where $\mathbf{A}$ is a dense matrix of order $N$. It will require $\mathcal{O}(N^2)$ operation for computing $\mathbf{A}\mathbf{z}$. FMM can be used to accelerate the computing.

The core of FMM is the multipole and local expansions of the integrals, and the translations of the expansion coefficients. In the FMM for computing $(S\varphi)(\mathbf{x})$, the multipole expansion, local expansion, and moment-to-local (M2L) translation are based on the Graf's addition theorem for Hankel function $H_m^{(1)}$, whereas moment-to-moment (M2M) and local-to-local (L2L) translations are based on that for Bessel function $J_m$. Graf's addition theorem is \cite{17,18}:
\begin{equation}
\mathscr{B}_m(|\mathbf{x}-\mathbf{y}|)e^{\pm{\mathrm{i}m}\theta_{\mathbf{x}-\mathbf{y}}}=\sum_{n=-\infty}^{\infty}\mathscr{B}_{m+n}(|\mathbf{x}|)e^{\pm{\mathrm{i}(m+n)}
\theta_\mathbf{x}}J_n(|\mathbf{y}|)e^{\mp{\mathrm{i}n}\theta_\mathbf{y}},\qquad |\mathbf{y}|<|\mathbf{x}|,
\end{equation}
another form is:
\begin{equation}
\mathscr{B}_m(|\mathbf{x}+\mathbf{y}|)e^{\pm{\mathrm{i}m}\theta_{\mathbf{x}+\mathbf{y}}}=\sum_{n=-\infty}^{\infty}\mathscr{B}_{m-n}(|\mathbf{x}|)e^{\pm{\mathrm{i}(m-n)}
\theta_\mathbf{x}}J_n(|\mathbf{y}|)e^{\pm{\mathrm{i}n}\theta_\mathbf{y}},\qquad |\mathbf{y}|<|\mathbf{x}|,
\end{equation}
where $m\in\mathds{Z}$, $\mathscr{B}$ denotes $J,Y,H^{(1)},H^{(2)}$ or any linear combination of these functions. When $\mathscr{B}=J$, the restriction $|\mathbf{y}|<|\mathbf{x}|$ is unnecessary. We denote the truncation error of (1) by
$$R^{\mathscr{B}}_{m,p}(\mathbf{x},\mathbf{y})=\left(\sum_{n=p+1}^\infty+\sum_{n=-\infty}^{-p-1}\right)\mathscr{B}_{m+n}(|\mathbf{x}|)
e^{\pm{\mathrm{i}(m+n)}\theta_\mathbf{x}}J_n(|\mathbf{y}|)e^{\mp{\mathrm{i}n}\theta_\mathbf{y}},$$
where $p$ is the truncation number. Obviously, the truncation error of (2) is $R^{\mathscr{B}}_{m,p}(\mathbf{x},-\mathbf{y})$. For convenience, let $\mathcal{H}^{\pm}_n(\mathbf{x})$ and $\mathcal{J}^{\pm}_n(\mathbf{x})$ be defined by
$$\mathcal{H}^{\pm}_n(\mathbf{x})=H^{(1)}_n(|\mathbf{x}|)e^{\pm{\mathrm{i}n}\theta_{\mathbf{x}}},\qquad
\mathcal{J}^{\pm}_n(\mathbf{x})=J_n(|\mathbf{x}|)e^{\pm{\mathrm{i}n}\theta_{\mathbf{x}}}.$$
Thus, Graf's addition theorem can be written as the form of $\mathcal{H}^{\pm}_n$ and $\mathcal{J}^{\pm}_n$.

Let us recall the truncation errors of the expansions and translations used in the FMM. Suppose that $\mathbf{y}_c$ is an expansion point close to the field point $\mathbf{y}$, that is, $|\mathbf{y}-\mathbf{y}_c|<|\mathbf{x}-\mathbf{y}_c|$. From (1), we have the following multipole expansion:
\begin{eqnarray}
\int_{\Delta\partial\Omega}H^{(1)}_0(k|\mathbf{x}-\mathbf{y}|)\varphi(\mathbf{y})\mathrm{d}s(\mathbf{y})\hspace{-0.6cm}&&=\int_{\Delta\partial\Omega}\Bigg(\sum_{n=-\infty}^{\infty}\mathcal{H}^+_n(k(\mathbf{x}-\mathbf{y}_c))
\mathcal{J}^-_n(k(\mathbf{y}-\mathbf{y}_c))\Bigg)\varphi(\mathbf{y})\mathrm{d}s(\mathbf{y})\nonumber\\
&&=\sum_{n=-p}^p\mathcal{H}^+_n(k(\mathbf{x}-\mathbf{y}_c))\mathbf{M}_{n}(\mathbf{y}_c)+E_{ME},
\end{eqnarray}
where
\begin{equation}
\mathbf{M}_n(\mathbf{y}_c)=\int_{\Delta\partial\Omega}\mathcal{J}^{-}_n(k(\mathbf{y}-\mathbf{y}_c))\varphi(\mathbf{y})\mathrm{d}s(\mathbf{y}), \qquad n\in\mathbb{Z}.
\end{equation}
are the multipole moments centered at $\mathbf{y}_c$, and
\begin{equation}
E_{ME}=\int_{\Delta\partial\Omega}R_{0,p}^{H}\big(k(\mathbf{x}-\mathbf{y}_c),k(\mathbf{y}-\mathbf{y}_c)\big)\varphi(\mathbf{y})\mathrm{d}s(\mathbf{y})
\end{equation}
is the truncation error of the multipole expansion.

When the multipole expansion point is moved from $\mathbf{y}_{c'}$ to $\mathbf{y}_c$, from (2), we have the following M2M translation:
\begin{eqnarray}
\mathbf{M}_n(\mathbf{y}_c)\hspace{-0.6cm}&&=\int_{\Delta\partial\Omega}\Bigg(\sum_{m=-\infty}^{\infty}\mathcal{J}^{-}_m(k(\mathbf{y}-\mathbf{y}_{c'}))
\mathcal{J}^{-}_{n-m}(k(\mathbf{y}_{c'}-\mathbf{y}_c))\Bigg)\varphi(\mathbf{y})\mathrm{d}s(\mathbf{y})\nonumber\\
&&=\sum_{m=-p}^{p}\mathcal{J}^{-}_{n-m}(k(\mathbf{y}_{c'}-\mathbf{y}_c))\mathbf{M}_m(\mathbf{y}_{c'})+EM_n,
\end{eqnarray}
where
\begin{equation}
EM_n=\int_{\Delta\partial\Omega}R_{n,p}^{J}\big(k(\mathbf{y}_{c'}-\mathbf{y}_c),-k(\mathbf{y}-\mathbf{y}_{c'})\big)\varphi(\mathbf{y})\mathrm{d}s(\mathbf{y})
\end{equation}
is the error of $\mathbf{M}_n(\mathbf{y}_c)$.

For the main part of (3), when $|\mathbf{x}-\mathbf{x}_c|<|\mathbf{x}_c-\mathbf{y}_c|$, we have the following local expansion:
\begin{eqnarray}
\sum_{n=-p}^p\mathcal{H}^+_n(k(\mathbf{x}-\mathbf{y}_c))\mathbf{M}_{n}(\mathbf{y}_c)\hspace{-0.6cm}&&=\sum_{n=-p}^{p}\Bigg(\sum_{m=-\infty}^{\infty}\mathcal{H}^+_{n-m}(k(\mathbf{x}_c-\mathbf{y}_c))
\mathcal{J}^+_m(k(\mathbf{x}-\mathbf{x}_c))\Bigg)\mathbf{M}_n(\mathbf{y}_c)\nonumber\\
&&=\sum_{m=-p}^p\mathcal{J}^+_m(k(\mathbf{x}-\mathbf{x}_c))\mathbf{L}_m(\mathbf{x}_c)+E_{ML},
\end{eqnarray}
where
\begin{equation}
\mathbf{L}_m(\mathbf{x}_c)=\sum_{n=-p}^{p}\mathcal{H}^+_{n-m}(k(\mathbf{x}_c-\mathbf{y}_c))\mathbf{M}_n(\mathbf{y}_c)
\end{equation}
is the M2L translation, $\mathbf{x}_c$ is the local expansion center, and
\begin{equation}
E_{ML}=\sum_{n=-p}^{p}R_{n,p}^{H}\big(k(\mathbf{x}_c-\mathbf{y}_c),-k(\mathbf{x}-\mathbf{x}_c)\big)\mathbf{M}_n(\mathbf{y}_c)
\end{equation}
is the error of the M2L translation.

The local expansion point can be moved from $\mathbf{x}_c$ to $\mathbf{x}_{c'}$, by (2), we have
\begin{eqnarray}
\sum_{m=-p}^p\mathcal{J}^+_m(k(\mathbf{x}-\mathbf{x}_c))\mathbf{L}_m(\mathbf{x}_c)\hspace{-0.6cm}&&=\sum_{m=-p}^{p}\left(\sum_{n=-\infty}^{\infty}\mathcal{J}^+_{m-n}(k(\mathbf{x}_{c'}-\mathbf{x}_c))
\mathcal{J}^+_n(k(\mathbf{x}-\mathbf{x}_{c'}))\right)\mathbf{L}_m(\mathbf{x}_c)\nonumber\\
&&=\sum_{n=-p}^{p}\mathcal{J}^+_n(k(\mathbf{x}-\mathbf{x}_{c'}))\mathbf{L}_n(\mathbf{x}_{c'})+E_{LL},
\end{eqnarray}
where
$$\mathbf{L}_n(\mathbf{x}_{c'})=\sum_{m=-p}^{p}\mathcal{J}^+_{n-m}(k(\mathbf{x}_{c'}-\mathbf{x}_c))\mathbf{L}_m(\mathbf{x}_c)$$
is the L2L translation and
\begin{equation}
E_{LL}=\sum_{m=-p}^{p}R_{m,p}^{J}\big(k(\mathbf{x}_{c'}-\mathbf{x}_c),-k(\mathbf{x}-\mathbf{x}_{c'})\big)\mathbf{L}_m(\mathbf{x}_c)
\end{equation}
is the error of the L2L translation.

From (3), (8) and (11), if $\mathbf{M}_n(\mathbf{y}_c)$ is calculated by (4) directly, then we have
\begin{equation}
\int_{\Delta\partial\Omega}H^{(1)}_0(k|\mathbf{x}-\mathbf{y}|)\varphi(\mathbf{y})\mathrm{d}s(\mathbf{y})=\sum_{n=-p}^{p}\mathcal{J}^+_n(k(\mathbf{x}-\mathbf{x}_{c'}))\mathbf{L}_n(\mathbf{x}_{c'})+E_{ME}+E_{ML}+E_{LL}.
\end{equation}
This implies that the error of the FMM is the sum of multipole expansion error, M2L error and L2L error. If $\mathbf{M}_n(\mathbf{y}_c)$ is obtained by the M2M translation (6), then we let
\begin{equation}
\mathbf{M}_n(\mathbf{y}_c)=\widetilde{\mathbf{M}}_n(\mathbf{y}_c,p)+EM_n,
\end{equation}
where $\widetilde{\mathbf{M}}_n(\mathbf{y}_c,p)$ is the approximation of $\mathbf{M}_n(\mathbf{y}_c)$ (i.e., the main part of (6)). Substituting (14) into the multipole expansion (3), an error caused by the M2M is generated, which is written as
\begin{equation}
E_{MM}=\sum_{n=-p}^p\mathcal{H}^+_n(k(\mathbf{x}-\mathbf{y}_c))EM_n.
\end{equation}
In addition, the multipole moment $\mathbf{M}_n(\mathbf{y}_c)$ in (9) and (10) is replaced by $\widetilde{\mathbf{M}}_n(\mathbf{y}_c,p)$. Thus, we have
\begin{equation}
\int_{\Delta\partial\Omega}H^{(1)}_0(k|\mathbf{x}-\mathbf{y}|)\varphi(\mathbf{y})\mathrm{d}s(\mathbf{y})=\sum_{n=-p}^{p}\mathcal{J}^+_n(k(\mathbf{x}-\mathbf{x}_{c'}))\mathbf{L}_n(\mathbf{x}_{c'})+E_{ME}+E_{MM}+E_{ML}+E_{LL}.
\end{equation}
From (5), (7), (10) and (12), we see that formulas of $E_{ME}$ and $E_{ML}$ include the truncation error $R^H_{m,p}$, while $E_{MM}$ and $E_{LL}$ include $R^J_{m,p}$. In Section 3, we will give sharp bounds on $R^J_{m,p}$ and $R^H_{m,p}$, and further estimate the above errors in Section 4.

\section{Bounds on $R^J_{m,p}(\mathbf{x},\mathbf{y})$ and $R^H_{m,p}(\mathbf{x},\mathbf{y})$}
In this section, we will estimate the bound on $R^{\mathscr{B}}_{m,p}(\mathbf{x},\mathbf{y})$ by the limiting forms of $J_n(z)$ and $Y_n(z)$. From $\mathscr{B}_{-n}=(-1)^n\mathscr{B}_n$, we have
\begin{equation}
\left|R^{\mathscr{B}}_{m,p}(\mathbf{x},\mathbf{y})\right|\leq\sum_{n=p+1}^\infty[|\mathscr{B}_{n+m}(|\mathbf{x}|)|+|\mathscr{B}_{n-m}(|\mathbf{x}|)|]|J_n(|\mathbf{y}|)|.
\end{equation}
This inequality shows that $R^{\mathscr{B}}_{-m,p}(\mathbf{x},\mathbf{y})$ and $R^{\mathscr{B}}_{m,p}(\mathbf{x},\mathbf{y})$ have the same upper bound. Hence, we will only consider the case $m\geq0$.

From Ref. \cite{18}, for all $n\in\mathds{N}$ and real number $z\geq0$, the following upper bounds on $|J_n(z)|$ hold.
\begin{eqnarray}
&&|J_n(z)|\leq\frac{1}{\Gamma(n+1)}\left(\frac{z}{2}\right)^n,\\[1ex]
&&|J_n(z)|\leq\left\{
 \begin{array}{ll}
 1, & n=0,\\[0.5ex]
 \displaystyle{\frac{1}{\sqrt{2}}}, & n\geq1.
 \end{array}\right.
\end{eqnarray}
We first give the bound on $R^J_{m,p}(\mathbf{x},\mathbf{y})$ in the following theorem.

{\theorem Let $m,p\in\mathds{N}$ and $\mathbf{x},\mathbf{y}\in\mathds{R}^2$. When $p\geq|\mathbf{y}|$,}
$$\left|R^J_{m,p}(\mathbf{x},\mathbf{y})\right|\leq\frac{4}{\Gamma(p+2)}\left(\frac{|\mathbf{y}|}{2}\right)^{p+1}.$$
\begin{proof} From (17), (18) and (19), we have
\begin{eqnarray*}
\left|R^J_{m,p}(\mathbf{x},\mathbf{y})\right|\hspace{-0.6cm}&&\leq\sum_{n=p+1}^\infty[|J_{n+m}(|\mathbf{x}|)|+|J_{n-m}(|\mathbf{x}|)|]|J_n(|\mathbf{y}|)|\\
&&\leq2\sum_{n=p+1}^\infty|J_n(|\mathbf{y}|)|\leq2\sum_{n=p+1}^\infty\frac{1}{\Gamma(n+1)}\left(\frac{|\mathbf{y}|}{2}\right)^n.
\end{eqnarray*}
Since when $n\geq|\mathbf{y}|$,
$$\frac{1}{\Gamma(n+2)}\left(\frac{|\mathbf{y}|}{2}\right)^{n+1}=\frac{|\mathbf{y}|}{2(n+1)}\frac{1}{\Gamma(n+1)}\left(\frac{|\mathbf{y}|}{2}\right)^n\leq\frac{1}{2}\frac{1}{\Gamma(n+1)}\left(\frac{|\mathbf{y}|}{2}\right)^n,$$
it follows that
$$\sum_{n=p+1}^\infty\frac{1}{\Gamma(n+1)}\left(\frac{|\mathbf{y}|}{2}\right)^n\leq\frac{1}{\Gamma(p+2)}\left(\frac{|\mathbf{y}|}{2}\right)^{p+1}\sum_{n=p+1}^\infty\frac{1}{2^{n-p-1}}
=\frac{2}{\Gamma(p+2)}\left(\frac{|\mathbf{y}|}{2}\right)^{p+1},$$
which proves the theorem.
\end{proof}

We now estimate the bound on $R^H_{m,p}(\mathbf{x},\mathbf{y})$. Since $H^{(1)}_n(z)=J_n(z)+\mathrm{i}Y_n(z)$, we should first consider the bound on $|Y_n(z)|$. For each integer $n>0$, as $z\rightarrow0$, the limiting form of $Y_n(z)$ \cite{17,18} is
$$Y_n(z)\sim-\left(\frac{2}{z}\right)^n\frac{\Gamma(n)}{\pi}.$$
We might as well define the function $C_n(z)$ by
\begin{equation}
C_n(z)=-Y_n(z)\left(\frac{z}{2}\right)^n\frac{\pi}{\Gamma(n)}.
\end{equation}
It is obvious that $C_n(z)\rightarrow1$ when $z\rightarrow0$. In addition, $C_n(z)$ also has the following properties.

{\lemma Suppose $n\in\mathds{N}, z\in\mathds{R}$ with $z>0$ and $n\geq{z}+1$. For fixed $n$, $C_n(z)$ is a strictly increasing function of $z$. For fixed $z$, as $n\rightarrow\infty$, $C_n(z)\rightarrow1$ and $1<C_{n+1}(z)<C_n(z)$.}
\begin{proof} Take the derivative of $C_n(z)$, we obtain
\begin{equation}
Y'_n(z)+\frac{n}{z}Y_n(z)=-C'_n(z)\left(\frac{2}{z}\right)^n\frac{\Gamma(n)}{\pi}.
\end{equation}
The recurrence relations of Bessel functions show that \cite{18}
\begin{equation}
Y_{n-1}(z)+Y_{n+1}(z)=\frac{2n}{z}Y_{n}(z),
\end{equation}
\begin{equation}
Y_{n-1}(z)-Y_{n+1}(z)=2Y'_{n}(z).
\end{equation}
Adding (22) and (23), by (21) we have
$$Y_{n-1}(z)=-C'_n(z)\left(\frac{2}{z}\right)^n\frac{\Gamma(n)}{\pi}.$$

Let $y_{n,1}$ be the first positive zero of $Y_n(z)$. From Ref. \cite{18}, $n<y_{n,1}$ and $Y_n(z)<0(0<z<y_{n,1})$. It follows that $Y_n(z)<0$ when $0<z\leq{n}$. Hence, when $0<z\leq{n}-1$, $Y_{n-1}(z)<0$, and then
$C'_n(z)>0$. This implies that $C_n(z)$ is a strictly increasing function. In addition, from $C_n(z)\rightarrow1(z\rightarrow0)$, we see that $C_n(z)>1$ when $0<z\leq{n}-1$.

Next, by the definition of $C_n(z)$ and (22), we have
$$C_{n+1}(z)-C_n(z)=\frac{\pi}{\Gamma(n+1)}\left(\frac{z}{2}\right)^{n+1}\left(-Y_{n+1}(z)+\frac{2n}{z}Y_{n}(z)\right)=\frac{\pi}{\Gamma(n+1)}\left(\frac{z}{2}\right)^{n+1}Y_{n-1}(z).$$
Since $Y_{n-1}(z)<0$ when $0<z\leq{n}-1$, it follows that $C_{n+1}(z)<C_n(z)$.

Finally, we consider the asymptotic behavior of $C_n(z)$. When $n\rightarrow\infty$, the asymptotic forms \cite{18}
$$\Gamma(n)\sim\sqrt{\frac{2\pi}{n}}\left(\frac{n}{e}\right)^n,\qquad Y_n(z)\sim{-\sqrt{\frac{2}{\pi{n}}}}\left(\frac{2n}{ez}\right)^n$$
show that
$$C_n(z)=-Y_n(z)\left(\frac{z}{2}\right)^n\frac{\pi}{\Gamma(n)}\sim1.$$
The proof is completed.
\end{proof}

From Lemma 1 and (19), we see that $|J_n(z)|\leq1<|Y_n(z)|$ when $n\geq{z}+1$, and moreover,
\begin{equation}
\big|H^{(1)}_n(z)\big|=\sqrt{J^2_n(z)+Y^2_n(z)}\leq\sqrt{2}|Y_n(z)|=\sqrt{2}C_n(z)\left(\frac{2}{z}\right)^n\frac{\Gamma(n)}{\pi}.
\end{equation}
Thus, the bound on $|H^{(1)}_n(z)|$ can be derived by the monotonicity of $C_n(z)$. It should be noted that Lemma 1 and the inequalities (18), (19), (24) also hold when $n$ is a positive real number.

We next give a novel estimate of the remainder term of the convergent power series:
\begin{equation}
\sum_{n=0}^{\infty}(n+t)^ar^n,\qquad 0<r<1,0\leq t<1,
\end{equation}
which will be used in the subsequent proof. In Ref. \cite{19}, an estimate about the remainder term for $t=0$ is proposed. In the following lemma, a simpler and sharper estimate is derived.

{\lemma Suppose $p\in\mathds{N}$ and $a\geq0$, $0<r<1$, $0\leq t<1$. When $p+t>-a/\ln{r}$,}
$$\sum_{n=p+1}^{\infty}(n+t)^ar^n\leq\frac{(p+t)^{a+1}r^{p}}{-(p+t)\ln{r}-a}.$$
\begin{proof}
Since when $x>-a/\ln{r}$, $x^ar^x$ is a strictly decreasing function of $x$. It follows that
\begin{equation}
\sum_{n=p+1}^{\infty}(n+t)^ar^n=r^{-t}\sum_{n=p+1}^{\infty}(n+t)^ar^{n+t}\leq r^{-t}\int_{p+t}^{+\infty}x^ar^x\mathrm{d}x
\end{equation}
when $p+t>-a/\ln{r}$, and furthermore,
$$\int_{p+t}^{+\infty}x^ar^x\mathrm{d}x=\int_{p+t}^{+\infty}x^a{e}^{x\ln{r}}\mathrm{d}x=\frac{1}{(-\ln{r})^{a+1}}\int_{-(p+t)\ln{r}}^{+\infty}y^ae^{-y}\mathrm{d}y=\frac{\Gamma\big(a+1,-(p+t)\ln{r}\big)}{(-\ln{r})^{a+1}},$$
where $\Gamma(\cdot,\cdot)$ is incomplete gamma function.

The asymptotic expansion of incomplete gamma function \cite{18} gives
$$\Gamma(a+1,z)=z^ae^{-z}\left(\sum_{i=0}^{m-1}\frac{a(a-1)\cdots(a-i+1)}{z^i}+\varepsilon_m(a+1,z)\right),$$
when $z>0$ and $m\geq{a}$,
$$\varepsilon_m(a+1,z)\leq\frac{|a(a-1)\cdots(a-m+1)|}{z^m}.$$
Let $m=[a]+1$, when $z>a$, we have
$$\Gamma(a+1,z)\leq{z}^ae^{-z}\sum_{i=0}^{[a]+1}\frac{a(a-1)\cdots(a-i+1)}{z^i}\leq{z}^ae^{-z}\sum_{i=0}^{[a]+1}\left(\frac{a}{z}\right)^i\leq\frac{{z}^{a+1}e^{-z}}{z-a}.$$
Thus, when $-(p+t)\ln{r}>a$,
\begin{equation}
\int_{p+t}^{+\infty}x^ar^x\mathrm{d}x=\frac{\Gamma\big(a+1,-(p+t)\ln{r}\big)}{(-\ln{r})^{a+1}}\leq\frac{(p+t)^{a+1}e^{(p+t)\ln{r}}}{-(p+t)\ln{r}-a}.
\end{equation}
(26) and (27) prove the lemma.
\end{proof}

The monotonicity of $|H_n(z)|$ given in the following lemma will be used to prove the subsequent theorem.

{\lemma Let $n\in\mathds{Z}$ and $z>0$. For fixed $n$, $|H_n(z)|$ is a strictly decreasing function of $z$. For fixed $z$, $|H_n(z)|$ is strictly increasing with the increase of $|n|$.}
\begin{proof} The conclusions for $n\geq0$ were proven in Ref. \cite{5}. In addition, since $|H_{-n}(z)|=|H_n(z)|$, it follows that the conclusions are also true for $n<0$. \end{proof}

On the basis of the above lemmas, we now give the bound on $R^H_{m,p}(\mathbf{x},\mathbf{y})$ in the following theorem.

{\theorem Let $m,p\in\mathds{N}$ and $\mathbf{x},\mathbf{y}\in\mathds{R}^2$ with $|\mathbf{y}|<|\mathbf{x}|$. When $p+m\geq|\mathbf{x}|$ and $p+m/2>-m/\ln{r}$,
$$\left|R^H_{m,p}(\mathbf{x},\mathbf{y})\right|\leq\frac{4\sqrt{2}C_{p+m+1}(|\mathbf{x}|)(2p+m)^mr^p}{\pi{|\mathbf{x}|}^m(-(2p+m)\ln{r}-2m)},$$
where $r=|\mathbf{y}|/|\mathbf{x}|$ and the function $C_n(\cdot)$ is defined by $(20)$.}
\begin{proof} From (17), (18), (24), Lemmas 1 and 3, when $p+m\geq|\mathbf{x}|$,
\begin{eqnarray}
\left|R^H_{m,p}(\mathbf{x},\mathbf{y})\right|\hspace{-0.6cm}&&\leq\sum_{n=p+1}^\infty[|H_{n+m}(|\mathbf{x}|)|+|H_{n-m}(|\mathbf{x}|)|]|J_n(|\mathbf{y}|)|\nonumber\\
&&\leq2\sum_{n=p+1}^\infty|H_{n+m}(|\mathbf{x}|)||J_n(|\mathbf{y}|)|\nonumber\\
&&\leq\frac{2^{m+\frac{3}{2}}}{\pi{|\mathbf{x}|}^m}\sum_{n=p+1}^\infty{C}_{n+m}(|\mathbf{x}|)\frac{\Gamma(n+m)}{\Gamma(n+1)}r^n\nonumber\\
&&\leq\frac{2^{m+\frac{3}{2}}C_{p+m+1}(|\mathbf{x}|)}{\pi{|\mathbf{x}|}^m}\sum_{n=p+1}^\infty\frac{\Gamma(n+m)}{\Gamma(n+1)}r^n,
\end{eqnarray}
where $r=|\mathbf{y}|/|\mathbf{x}|<1$. The inequality of arithmetic and geometric means gives
$$\frac{\Gamma(n+m)}{\Gamma(n+1)}=(n+m-1)(n+m-2)\cdots(n+1)\leq\left(n+\frac{m}{2}\right)^{m-1}.$$
Let $m/2=l+t$, where $l$ is an integer and $t=0$ or $1/2$. By Lemma 2, when $p+l+t>-m/\ln{r}$,
\begin{eqnarray}
\sum_{n=p+1}^\infty\frac{\Gamma(n+m)}{\Gamma(n+1)}r^n\hspace{-0.6cm}&&\leq\sum_{n=p+1}^\infty\left(n+l+t\right)^{m-1}r^n\nonumber\\
&&=r^{-l}\sum_{n=p+l+1}^\infty\left(n+t\right)^{m-1}r^n\nonumber\\
&&\leq\frac{r^{-l}}{p+l+t+1}\sum_{n=p+l+1}^\infty\left(n+t\right)^mr^n\nonumber\\
&&\leq\frac{(p+l+t)^mr^p}{-(p+l+t)\ln{r}-m}.
\end{eqnarray}
From (28) and (29), we obtain
$$\left|R^H_{m,p}(\mathbf{x},\mathbf{y})\right|\leq\frac{4\sqrt{2}C_{p+m+1}(|\mathbf{x}|)(2p+m)^mr^p}{\pi{|\mathbf{x}|}^m(-(2p+m)\ln{r}-2m)},$$
which proves the theorem.
\end{proof}

At the end of this section, we perform some numerical experiments to test the bound given in the above theorem. We compare the bound derived here with those from Ref. \cite{19}.

In Fig.1, $|R^H_{m,p}(\mathbf{x},\mathbf{y})|$ and its bound are plotted as functions of $p$ and $m$ respectively. We see that the bound given here is in close agreement with the exact value, for all $m$ and $p$. In addition, the bound is sharper than the previous result, especially for smaller $p$ and larger $m$.
\begin{figure}[ht]
\centering
\scalebox{0.55}{\includegraphics{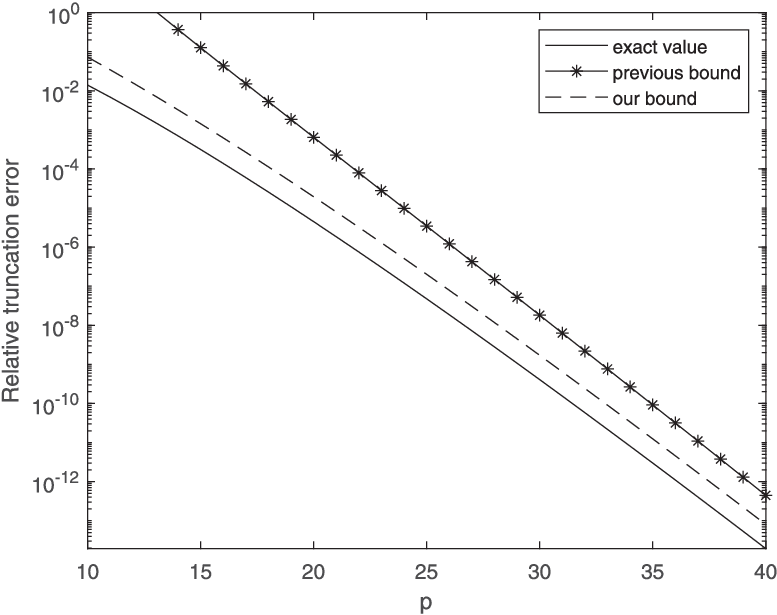}}\qquad\scalebox{0.55}{\includegraphics{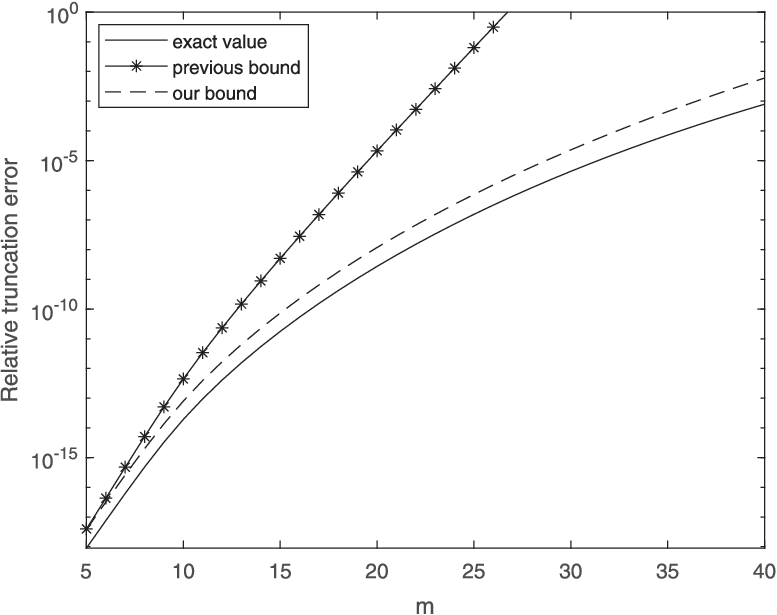}}\\
{\small Fig.1. $|R^H_{10,p}(\mathbf{x},\mathbf{y})|$ (left) and $|R^H_{m,40}(\mathbf{x},\mathbf{y})|$ (right), $|\mathbf{x}|=10,|\mathbf{y}|=3,\theta_\mathbf{x}=\pi/3,\theta_\mathbf{y}=\pi/6$.}
\end{figure}

Theorem 2 shows that the bound on $R^H_{m,p}(\mathbf{x},\mathbf{y})$ is increasing with the increase of $m$. From this and the error of the M2L translation (10), we see that the estimate of $R^H_{p,p}(\mathbf{x},\mathbf{y})$ is crucial to the FMM error. In Fig.2, $R^H_{p,p}(\mathbf{x},\mathbf{y})$ and its bound are plotted as functions of $p$. It is shown that the bound given here is in close agreement with the exact value.
\begin{figure}[ht]
\centering
\scalebox{0.6}{\includegraphics{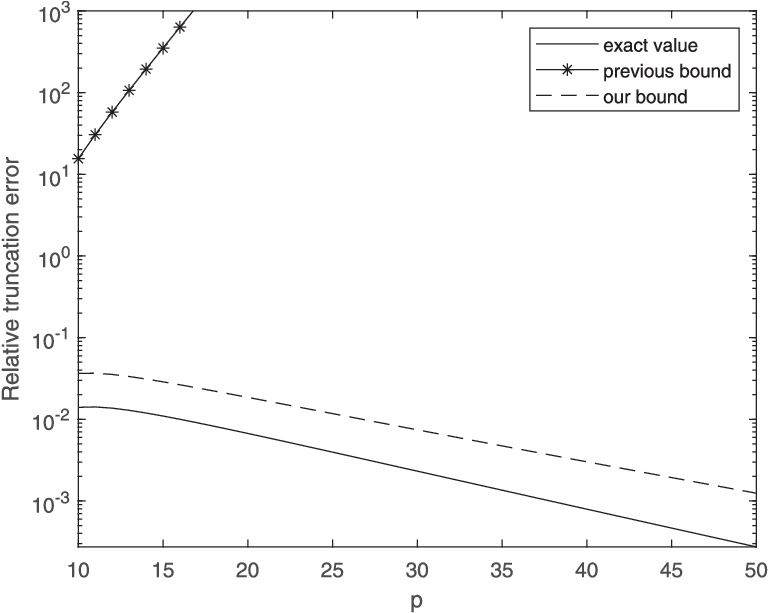}}\\
{\small Fig.2. $|R^H_{p,p}(\mathbf{x},\mathbf{y})|$ and its bound, $|\mathbf{x}|=10,|\mathbf{y}|=3,\theta_\mathbf{x}=\pi/3,\theta_\mathbf{y}=\pi/6$.}
\end{figure}

As can be seen from the above numerical experiments, the bound given here is very sharp in all cases. This is really because we use the limiting forms of the Bessel functions and give a very sharp bound for the remainder term of the series (25).

\section{Error bound of the FMM }
In this section, the FMM error will be estimated by using the bounds on $R^J_{m,p}(\mathbf{x},\mathbf{y})$ and $R^H_{m,p}(\mathbf{x},\mathbf{y})$ proposed in the previous section. And then, the result will be applied to the adaptive tree structure, the specific convergence order of the FMM error is obtained.

\subsection{Error bound of the FMM for $(S\varphi)(\mathbf{x})$}
We now give the estimates of the bounds on $E_{ME},E_{MM},E_{ML}$ and $E_{LL}$ mentioned in Section 2. Their detailed formulas are as follows:
\begin{eqnarray*}
&&E_{ME}=\int_{\Delta\partial\Omega}R_{0,p}^{H}\big(k(\mathbf{x}-\mathbf{y}_c),k(\mathbf{y}-\mathbf{y}_c)\big)\varphi(\mathbf{y})\mathrm{d}s(\mathbf{y}),\\[1ex]
&&E_{MM}=\sum_{n=-p}^p\mathcal{H}^+_n(k(\mathbf{x}-\mathbf{y}_c))EM_n,\\
&&E_{ML}=\sum_{n=-p}^{p}R_{n,p}^{H}\big(k(\mathbf{x}_c-\mathbf{y}_c),-k(\mathbf{x}-\mathbf{x}_c)\big)\mathbf{M}_n(\mathbf{y}_c),\\
&&E_{LL}=\sum_{n=-p}^{p}R_{n,p}^{J}\big(k(\mathbf{x}_{c'}-\mathbf{x}_c),-k(\mathbf{x}-\mathbf{x}_{c'})\big)\mathbf{L}_n(\mathbf{x}_c).
\end{eqnarray*}
Without loss of generality, we assume that $\varphi(\mathbf{y})$ is a continuous function on $\Delta\partial\Omega$, and $\varphi(\mathbf{y})\not\equiv0$.

Before estimating the above errors, we first give the bound on $EM_n$. From (7) and Theorem 1, when $p\geq k|\mathbf{y}-\mathbf{y}_{c'}|$,
\begin{eqnarray}
\big|EM_n\big|\hspace{-0.6cm}&&\leq\int_{\Delta\partial\Omega}\left|R_{n,p}^{J}\big(k(\mathbf{y}_{c'}-\mathbf{y}_c),-k(\mathbf{y}-\mathbf{y}_{c'})\big)\right||\varphi(\mathbf{y})|\mathrm{d}s(\mathbf{y})\nonumber\\
&&\leq\frac{4}{\Gamma(p+2)}\left(\frac{kt}{2}\right)^{p+1}\int_{\Delta\partial\Omega}|\varphi(\mathbf{y})|\mathrm{d}s(\mathbf{y})\nonumber\\
&&\leq\frac{A(kt)^{p+1}}{2^{p-1}\Gamma(p+2)},
\end{eqnarray}
where $t=\sup\{|\mathbf{y}-\mathbf{y}_{c'}|\}$, $A=\max\{|\varphi(\mathbf{y})||\Delta\partial\Omega|\}$ and $|\Delta\partial\Omega|$ denotes the length of $\Delta\partial\Omega$.

{\theorem Let $\mathbf{x},\mathbf{y},\mathbf{x}_c,\mathbf{y}_c,\mathbf{x}_{c'},\mathbf{y}_{c'}$ satisfy the conditions of the expansions and translations of the FMM. Then, the bounds on truncation errors of the FMM for $(S\varphi)(\mathbf{x})$ are given by
\begin{eqnarray*}
&&|E_{ME}|\leq\frac{2\sqrt{2}AC_{p+1}(k|\mathbf{x}-\mathbf{y}_c|)\mu^p}{-\pi{p}\ln{\mu}}=\mathcal{O}\left(\frac{\mu^p}{p}\right),\\
&&|E_{MM}|\leq\frac{4\sqrt{2}AC_p(k|\mathbf{x}-\mathbf{y}_c|)kt\eta^p}{\pi{p}}=\mathcal{O}\left(\frac{\eta^p}{p}\right),\\
&&|E_{ML}|\leq\frac{12\sqrt{2}AC_{p+1}(k|\mathbf{x}_c-\mathbf{y}_c|)\lambda^p}{-\pi{p}(3\ln\gamma+2)}=\mathcal{O}\left(\frac{\lambda^p}{p}\right),\\
&&|E_{LL}|\leq\frac{8\sqrt{2}AC_p(k\tau)k|\mathbf{x}-\mathbf{x}_{c'}|\nu^p}{\pi{p}}=\mathcal{O}\left(\frac{\nu^p}{p}\right),
\end{eqnarray*}
for $\gamma<e^{-2/3}$ and $p\geq\max\{k|\mathbf{x}-\mathbf{y}_c|+1,k|\mathbf{x}_c-\mathbf{y}_c|,k|\mathbf{x}_c-\mathbf{y}|+1\}$, where $A=\max\{|\varphi(\mathbf{y})||\Delta\partial\Omega|\}$, $\mu=\sup\{|\mathbf{y}-\mathbf{y}_c|/|\mathbf{x}-\mathbf{y}_c|\}$, $\eta=\sup\{|\mathbf{y}-\mathbf{y}_{c'}|/|\mathbf{x}-\mathbf{y}_c|\}$, $t=\sup\{|\mathbf{y}-\mathbf{y}_{c'}|\}$, $\rho=\sup\{|\mathbf{y}-\mathbf{y}_c|/|\mathbf{x}_c-\mathbf{y}_c|\}$, $\gamma=|\mathbf{x}-\mathbf{x}_c|/|\mathbf{x}_c-\mathbf{y}_c|$, $\nu=\sup\left\{|\mathbf{x}-\mathbf{x}_{c'}|/|\mathbf{x}_c-\mathbf{y}|\right\}$, $\tau=\sup\{|\mathbf{x}_c-\mathbf{y}|\}$ and $\lambda=\gamma{e}^{\frac{3\rho}{2}}$.}
\begin{proof} (1). Let us first prove the result about $E_{ME}$. From Theorem 2, when $p\geq k|\mathbf{x}-\mathbf{y}_c|$,
\begin{eqnarray}
|E_{ME}|\hspace{-0.6cm}&&\leq\int_{\Delta\partial\Omega}\left|R_{0,p}^{H}\big(k(\mathbf{x}-\mathbf{y}_c),k(\mathbf{y}-\mathbf{y}_c)\big)\right||\varphi(\mathbf{y})|\mathrm{d}s(\mathbf{y})\nonumber\\[0.5ex]
&&\leq\frac{2\sqrt{2}C_{p+1}(k|\mathbf{x}-\mathbf{y}_c|)\mu^p}{-\pi{p}\ln{\mu}}\int_{\Delta\partial\Omega}|\varphi(\mathbf{y})|\mathrm{d}s(\mathbf{y})\nonumber\\
&&\leq\frac{2\sqrt{2}AC_{p+1}(k|\mathbf{x}-\mathbf{y}_c|)\mu^p}{-\pi{p}\ln{\mu}},
\end{eqnarray}
where $\mu=\sup\{|\mathbf{y}-\mathbf{y}_c|/|\mathbf{x}-\mathbf{y}_c|\}$ and $A=\max\{|\varphi(\mathbf{y})||\Delta\partial\Omega|\}$. In addition, from $C_{p+1}(z)\rightarrow1(p\rightarrow\infty)$, we prove the result about $E_{ME}$.

(2). Next for $E_{MM}$. From (24), (30) and Lemma 3, when $p\geq k|\mathbf{x}-\mathbf{y}_c|+1$, we have
\begin{eqnarray}
|E_{MM}|\hspace{-0.6cm}&&\leq\sum_{n=-p}^p\left|\mathcal{H}^+_n(k(\mathbf{x}-\mathbf{y}_c))\right||EM_n|\nonumber\\
&&\leq\left|H_{p}^{(1)}(k|\mathbf{x}-\mathbf{y}_c|)\right|\sum_{n=-p}^p|EM_n|\nonumber\\
&&\leq\frac{4\sqrt{2}AC_p(k|\mathbf{x}-\mathbf{y}_c|)kt\eta^p}{\pi{p}},
\end{eqnarray}
where $\eta=\sup\{|\mathbf{y}-\mathbf{y}_{c'}|/|\mathbf{x}-\mathbf{y}_c|\}$, which proves the second result.

(3). We prove the result for $E_{ML}$. Substituting the multipole moment (4) into $E_{ML}$, by Theorem 2 and (18), when $p\geq k|\mathbf{x}_c-\mathbf{y}_c|$, we have
\begin{eqnarray}
|E_{ML}|\hspace{-0.6cm}&&\leq\sum_{n=-p}^{p}\left|R_{n,p}^{H}\big(k(\mathbf{x}_c-\mathbf{y}_c),-k(\mathbf{x}-\mathbf{x}_c)\big)\right|\left|\mathbf{M}_n(\mathbf{y}_c)\right|\nonumber\\
&&\leq2\sum_{n=0}^{p}\frac{4\sqrt{2}C_{p+n+1}(k|\mathbf{x}_c-\mathbf{y}_c|)(2p+n)^n\gamma^p}{\pi(k|\mathbf{x}_c-\mathbf{y}_c|)^n(-(2p+n)\ln\gamma-2n)}
\int_{\Delta\partial\Omega}\frac{(k|\mathbf{y}-\mathbf{y}_c|)^n}{2^n\Gamma(n+1)}|\varphi(\mathbf{y})|\mathrm{d}s(\mathbf{y})
\end{eqnarray}
for $\ln\gamma<-2/3$, where $\gamma=|\mathbf{x}-\mathbf{x}_c|/|\mathbf{x}_c-\mathbf{y}_c|$. Lemma 1 shows that
$$C_{p+n+1}(k|\mathbf{x}_c-\mathbf{y}_c|)\leq C_{p+1}(k|\mathbf{x}_c-\mathbf{y}_c|).$$
Thus, we have
\begin{eqnarray}
|E_{ML}|\hspace{-0.6cm}&&\leq\frac{8\sqrt{2}C_{p+1}(k|\mathbf{x}_c-\mathbf{y}_c|)\gamma^p}{\pi}\sum_{n=0}^{p}\frac{(2p+n)^{n-1}}{(-\ln\gamma-\frac{2n}{2p+n})\Gamma(n+1)}
\int_{\Delta\partial\Omega}\frac{(|\mathbf{y}-\mathbf{y}_c|)^n}{(2|\mathbf{x}_c-\mathbf{y}_c|)^n}|\varphi(\mathbf{y})|\mathrm{d}s(\mathbf{y})\nonumber\\
&&\leq\frac{24\sqrt{2}AC_{p+1}(k|\mathbf{x}_c-\mathbf{y}_c|)\gamma^p}{\pi(-3\ln\gamma-2)}\sum_{n=0}^{p}\frac{(2p+n)^{n-1}}{\Gamma(n+1)}\left(\frac{\rho}{2}\right)^n\nonumber\\
&&\leq\frac{12\sqrt{2}AC_{p+1}(k|\mathbf{x}_c-\mathbf{y}_c|)\gamma^p}{-\pi{p}(3\ln\gamma+2)}\sum_{n=0}^{p}\frac{1}{\Gamma(n+1)}\left(\frac{3p\rho}{2}\right)^n\nonumber\\
&&\leq\frac{12\sqrt{2}AC_{p+1}(k|\mathbf{x}_c-\mathbf{y}_c|)\gamma^pe^{\frac{3p\rho}{2}}}{-\pi{p}(3\ln\gamma+2)},
\end{eqnarray}
where $\rho=\sup\{|\mathbf{y}-\mathbf{y}_c|/|\mathbf{x}_c-\mathbf{y}_c|\}$. This proves the third result.

(4). Finally, let us prove the result about $E_{LL}$. The M2L translation (9), multipole moment (4) and Graf's addition theorem (1) give
\begin{eqnarray}
\mathbf{L}_n(\mathbf{x}_c)\hspace{-0.6cm}&&=\sum_{m=-p}^{p}\mathcal{H}^+_{m-n}(k(\mathbf{x}_c-\mathbf{y}_c))\mathbf{M}_m(\mathbf{y}_c)\nonumber\\
&&=\int_{\Delta\partial\Omega}\sum_{m=-p}^{p}\mathcal{H}^+_{m-n}(k(\mathbf{x}_c-\mathbf{y}_c))\mathcal{J}^{-}_m(k(\mathbf{y}-\mathbf{y}_c))\varphi(\mathbf{y})\mathrm{d}s(\mathbf{y})\nonumber\\
&&=\int_{\Delta\partial\Omega}\left[\mathcal{H}^+_{-n}(k(\mathbf{x}_c-\mathbf{y}))-R_{-n,p}^{H}\big(k(\mathbf{x}_c-\mathbf{y}_c),k(\mathbf{y}-\mathbf{y}_c)\big)\right]\varphi(\mathbf{y})\mathrm{d}s(\mathbf{y}).
\end{eqnarray}
Since $R_{-n,p}^{H}$ is the truncation error of the expansion about $\mathcal{H}^+_{-n}$, it follows that $|R_{-n,p}^{H}|<|\mathcal{H}^+_{-n}|$ when $p$ is large. Fig.1 shows that the inequality holds when $p\geq n$.
Thus, we let
\begin{equation}
\left|\mathbf{L}_n(\mathbf{x}_c)\right|\leq2\int_{\Delta\partial\Omega}\left|\mathcal{H}^+_{-n}(k(\mathbf{x}_c-\mathbf{y}))\right||\varphi(\mathbf{y})|\mathrm{d}s(\mathbf{y}).
\end{equation}
Now, from Theorem 1, Lemma 3 and (24), when $p\geq k|\mathbf{x}_c-\mathbf{y}|+1$,
\begin{eqnarray}
|E_{LL}|\hspace{-0.6cm}&&\leq\sum_{n=-p}^{p}\left|R_{n,p}^{J}\big(k(\mathbf{x}_{c'}-\mathbf{x}_c),-k(\mathbf{x}-\mathbf{x}_{c'})\big)\right|\left|\mathbf{L}_n(\mathbf{x}_c)\right|\nonumber\\
&&\leq\frac{8}{\Gamma(p+2)}\left(\frac{k|\mathbf{x}-\mathbf{x}_{c'}|}{2}\right)^{p+1}\sum_{n=-p}^{p}\int_{\Delta\partial\Omega}\left|H^{(1)}_n(k|\mathbf{x}_c-\mathbf{y}|)\right||\varphi(\mathbf{y})|\mathrm{d}s(\mathbf{y})\nonumber\\
&&\leq\frac{8(2p+1)}{\Gamma(p+2)}\left(\frac{k|\mathbf{x}-\mathbf{x}_{c'}|}{2}\right)^{p+1}\int_{\Delta\partial\Omega}\left|H^{(1)}_p(k|\mathbf{x}_c-\mathbf{y}|)\right||\varphi(\mathbf{y})|\mathrm{d}s(\mathbf{y})\nonumber\\
&&\leq\frac{4\sqrt{2}(2p+1)(k|\mathbf{x}-\mathbf{x}_{c'}|)^{p+1}}{\pi(p+1)p}\int_{\Delta\partial\Omega}\frac{C_p(k|\mathbf{x}_c-\mathbf{y}|)}{(k|\mathbf{x}_c-\mathbf{y}|)^p}|\varphi(\mathbf{y})|\mathrm{d}s(\mathbf{y})\nonumber\\
&&\leq\frac{8\sqrt{2}AC_p(k\tau)k|\mathbf{x}-\mathbf{x}_{c'}|\nu^p}{\pi{p}},
\end{eqnarray}
in which $\tau=\sup\{|\mathbf{x}_c-\mathbf{y}|\}$ and $\nu=\sup\left\{|\mathbf{x}-\mathbf{x}_{c'}|/|\mathbf{x}_c-\mathbf{y}|\right\}$. The proof is completed.
\end{proof}

The condition $\gamma<e^{-2/3}(\approx0.513)$ is necessary in Theorem 3. In the tree structure of the FMM algorithm, the value of $\gamma$ is less than $1/2$. Thus, the theorem is valid for the algorithm.

In fact, Christophe\cite{12} gave an estimate for the upper bound of $E_{ML}$, that is
\begin{equation}
\left|E_{ML}\right|\leq{C}p\left(\frac{2e\rho}{1-\gamma}\right)^{p-q},
\end{equation}
where $q$ is a nonnegative integer and $C$ is an unknown constant. The estimation is based on the conditions:
$$0<\gamma<\frac{1}{1+2e},\qquad 0<\rho<\frac{1}{1+2e}.$$
However, in the adaptive square tree structure, the maximum values of $\gamma$ and $\rho$ are $\sqrt{2}/4$ and $2/\sqrt{13}$ respectively (see next section). It follows that the estimate (38) is invalid in most cases. Theorem 3 of this paper gives
$$|E_{ML}|=\mathcal{O}\left(\frac{\lambda^p}{p}\right),$$
where $\lambda=\gamma{e}^{\frac{3\rho}{2}}$. By simple calculation, we see that $\lambda$ is much smaller than $2e\rho/(1-\gamma)$. Hence, our estimate is obviously novel, sharp and valid. This is because we give a very sharp estimate for $R_{n,p}^{H}$ in Theorem 2.

It should be noted that, if the multipole moment $\mathbf{M}_n(\mathbf{y}_c)$ is obtained by the M2M translation, then $\mathbf{M}_n(\mathbf{y}_c)=\widetilde{\mathbf{M}}_n(\mathbf{y}_c,p)+EM_n$. It follows that, in M2L and $E_{ML}$, $\mathbf{M}_n(\mathbf{y}_c)$ should be replaced by $\widetilde{\mathbf{M}}_n(\mathbf{y}_c,p)$. However, since $EM_n\rightarrow0$ (see (30)) and $\mathbf{M}_n(\mathbf{y}_c)$ is nonzero constant (when $k|\mathbf{y}-\mathbf{y}_c|$ happens to be the zero of $J_n(z)$, $\mathbf{M}_n(\mathbf{y}_c)=0$, but this is a small probability event and is not considered), we omitted the difference between $\mathbf{M}_n(\mathbf{y}_c)$ and $\widetilde{\mathbf{M}}_n(\mathbf{y}_c,p)$ in the proof of Theorem 3. In other words, some tedious and inconsequential parts are omitted.

\subsection{Application in the tree structure}
In this section, we will apply the estimates given in the above section to the tree structure and derive the convergence order of the FMM error. We have analyzed that the FMM error is the sum of $E_{ME}, E_{MM}, E_{ML}$ and $E_{LL}$, thus from Theorem 3, the convergence order of the FMM is determined by the maximum of $\mu, \eta, \lambda, \nu$.

In FMM, the expansions and translations are accomplished in the tree structure. In 2-D, the square quadtree structure is common used. There are two types of tree structures, called adaptive and nonadaptive. The so-called adaptive tree can automatically adjust its structure according to the element distributions. It is more efficient for BEM models with nonuniform element distributions. See Refs. \cite{7,8} for detailed introductions of the adaptive tree structure.

For illustration purposes, we denote the square by its center point. Suppose that $\mathbf{x}$ and $\mathbf{y}$ are contained in the squares centered on $\mathbf{x}_c$ and $\mathbf{y}_c$, respectively. In the nonadaptive tree structure, $\mathbf{x}_c$ and $\mathbf{y}_c$ are at the same layer of the tree, but it is not so in the adaptive one. This will affect the values of $\mu, \eta, \lambda$ and $\nu$.
\begin{figure}[htbp]
\centering
\scalebox{0.8}{\includegraphics{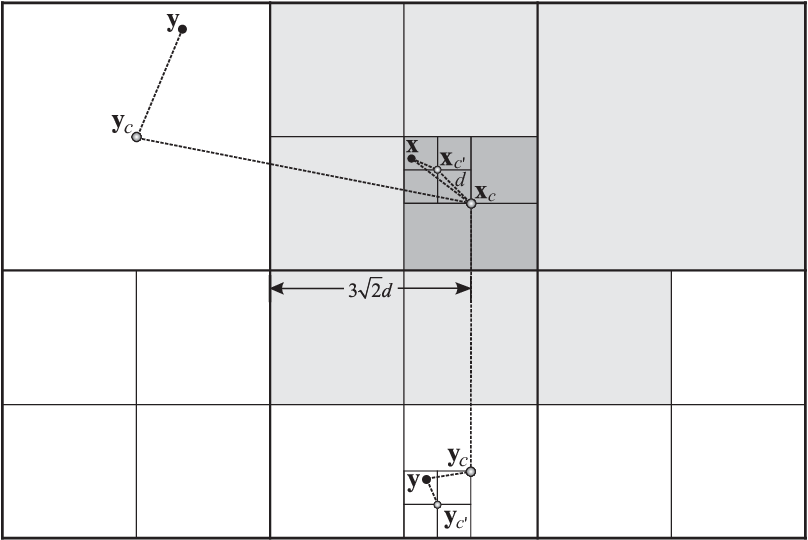}}\\
{\small Fig.3. The geometry of $\mathbf{x},\mathbf{y},\mathbf{x}_c,\mathbf{y}_c,\mathbf{x}_{c'},\mathbf{y}_{c'}$ in the adaptive tree structure.}
\end{figure}

Fig.3 shows the geometry of $\mathbf{x},\mathbf{y},\mathbf{x}_c,\mathbf{y}_c,\mathbf{x}_{c'},\mathbf{y}_{c'}$ in the adaptive tree structure, where $\mathbf{x}_{c'}$ and $\mathbf{y}_{c'}$ are the sons of $\mathbf{x}_c$ and $\mathbf{y}_c$, respectively. Let the layers of $\mathbf{x}_c$ and $\mathbf{y}_c$ be $L_{\mathbf{x}_c}$ and $L_{\mathbf{y}_c}$, respectively. Set $|\mathbf{x}_c-\mathbf{x}_{c'}|=d$, the we have $|\mathbf{x}-\mathbf{x}_{c'}|\leq d$, $|\mathbf{x}-\mathbf{x}_c|\leq2d$, $|\mathbf{x}_c-\mathbf{y}|\geq3\sqrt{2}d$ and
\begin{eqnarray*}
&&|\mathbf{x}-\mathbf{y}_c|\geq\left\{
 \begin{array}{ll}
 3\sqrt{2}d, & L_{\mathbf{y}_c}=L_{\mathbf{x}_c},\\[0.5ex]
 4\sqrt{2}d, & L_{\mathbf{y}_c}=L_{\mathbf{x}_c}-1.
 \end{array}\right. \qquad
|\mathbf{x}_c-\mathbf{y}_c|\geq\left\{
 \begin{array}{ll}
 4\sqrt{2}d, & L_{\mathbf{y}_c}=L_{\mathbf{x}_c},\\[0.5ex]
 2\sqrt{13}d, & L_{\mathbf{y}_c}=L_{\mathbf{x}_c}-1.
 \end{array}\right.\\
&&|\mathbf{y}-\mathbf{y}_c|\leq\left\{
 \begin{array}{ll}
 2d, & L_{\mathbf{y}_c}=L_{\mathbf{x}_c},\\[0.5ex]
 4d, & L_{\mathbf{y}_c}=L_{\mathbf{x}_c}-1.
 \end{array}\right.\qquad~~~~
|\mathbf{y}-\mathbf{y}_{c'}|\leq\left\{
\begin{array}{ll}
 d, & L_{\mathbf{y}_c}=L_{\mathbf{x}_c},\\[0.5ex]
 \backslash, & L_{\mathbf{y}_c}=L_{\mathbf{x}_c}-1.
 \end{array}\right.
\end{eqnarray*}
Note that if $L_{\mathbf{y}_c}=L_{\mathbf{x}_c}-1$, $\mathbf{y}_c$ is a leaf cell and has no son, then $\mathbf{y}_{c'}$ does not exist.

Substituting the above inequalities into $\mu,\eta,\gamma,\rho,\lambda,\nu$ mentioned in Theorem 3, we can derive their values as shown in Table 1.
\begin{center}
{Table 1. The values of $\mu,\eta,\gamma,\rho,\lambda,\nu$.}\\
{\begin{tabular}{|l|c|c|c|c|c|c|}\hline
                 & $\mu$ & $\eta$ & $\gamma$ & $\rho$ & $\lambda$ & $\nu$ \\\hline
 ~$L_{\mathbf{y}_c}=L_{\mathbf{x}_c}$ & ~~~$\sqrt{2}/3$~~~ & ~~~$\sqrt{2}/6$~~~ & ~~~$\sqrt{2}/4$~~~ & ~~~$\sqrt{2}/4$~~~ & $\sqrt{2}e^{\frac{3}{4\sqrt{2}}}/4$ & ~~~$\sqrt{2}/6$~~~  \\\hline
 ~$L_{\mathbf{y}_c}=L_{\mathbf{x}_c}-1$ & $\sqrt{2}/2$ & $\backslash$ & $1/\sqrt{13}$ & $2/\sqrt{13}$ & $e^{\frac{3}{\sqrt{13}}}/\sqrt{13}$ & ~$\sqrt{2}/6$~  \\\hline
\end{tabular}}
\end{center}

From Table 1, the convergence order of the FMM error with square tree structure is derived. See Table 2 for details.
\begin{center}
{Table 2. Convergence order of the FMM error.}\\
{\begin{tabular}{|l|l|c|c|c|}\hline
                & convergence order & $L_{\mathbf{y}_c}=L_{\mathbf{x}_c}$ & $L_{\mathbf{y}_c}=L_{\mathbf{x}_c}-1$ \\\hline
 ~~$E_{ME}$~ & ~~~~~~$p^{-1}\mu^p$~~ & ~$p^{-1}(0.4714)^p$~ & ~\cellcolor[gray]{0.8}{$p^{-1}(0.7071)^p$}~ \\\hline
 ~~$E_{MM}$~ & ~~~~~~$p^{-1}\eta^p$~ & ~$p^{-1}(0.2357)^p$~ & ~$\backslash$~ \\\hline
 ~~$E_{ML}$~ & ~~~~~~$p^{-1}\lambda^p$~ & ~\cellcolor[gray]{0.8}{$p^{-1}(0.6009)^p$}~ & ~$p^{-1}(0.6374)^p$~ \\\hline
 ~~$E_{LL}$~ & ~~~~~~$p^{-1}\nu^p$~ & ~$p^{-1}(0.2357)^p$~ & ~$p^{-1}(0.2357)^p$~ \\\hline
\end{tabular}}
\end{center}

The result shows that $E_{ME}$ and $E_{ML}$ are the two main parts of the FMM error. When $L_{\mathbf{y}_c}=L_{\mathbf{x}_c}$, $E_{ME}$ is smaller than $E_{ML}$, while the result is opposite when $L_{\mathbf{y}_c}=L_{\mathbf{x}_c}-1$. In addition, the FMM with nonadaptive tree structure has higher convergence order than that with adaptive one, although the latter is more efficient.

It should be noted that the convergence order shown in Table 2 is theoretical. In the real numerical examples, random distribution of $\mathbf{x}$ and $\mathbf{y}$ leads to the values of $\mu,\eta,\lambda,\nu$ which is smaller than those given in the table. Therefore, the convergence order of the FMM is determined by their specific values in the program.

We perform a numerical experiment to validate the above results. The boundary considered here is kite-shaped, with the parametric representation
$$\partial\Omega:(\cos{t}+0.65\cos2t-0.65,1.5\sin{t}).\qquad 0\leq{t}\leq2\pi.$$
The initial square is $[-1.8,1.8]\times[-1.8,1.8]$. We discretize the boundary into $N=1000$ constant elements, and then construct an adaptive and nonadaptive tree structures respectively. The leaf cell covers up to $\ln{N}$ points. In the adaptive tree, the values of $\mu$ and $\lambda$ are $0.6761$ and $0.5483$ respectively, while they are $0.4315$ and $0.5262$ in the nonadaptive tree.

In each tree structure, we choose a pair of source and field points that can represent the upper bound of the algorithm error. See Table 3 for details.
\begin{center}
{Table 3. Two pairs of representative points.}\\
{\begin{tabular}{|l|c|c|c|c|}\hline
                 & source point & field point & $\mu$ & $\lambda$ \\\hline
 adaptive & ~~~$\mathbf{x}_{97}$~~~ & ~~~~~$\mathbf{x}_{102}$~~~~~ & ~~~$0.6761$~~~ & ~~~$0.4812$~~~  \\\hline
 nonadaptive & $\mathbf{x}_{997}$ & $\mathbf{x}_{1000}$ & $0.4030$ & $0.5262$  \\\hline
\end{tabular}}
\end{center}
We will give the error in the above two cases. Let the wave number $k=5$ and the density function $\varphi(\mathbf{y})\equiv1$. An easy computation shows that
$$A=\max\{|\varphi(\mathbf{y})||\Delta\partial\Omega|\}\approx2.2718\times10^{-3}.$$

In Fig.4, for the above two cases, the FMM error and its bound are plotted as functions of $p$. It is seen that the proposed bound is valid for both adaptive and nonadaptive methods. And the nonadaptive method has higher convergence order than the adaptive one.
\begin{figure}[ht]
\centering
\scalebox{0.55}{\includegraphics{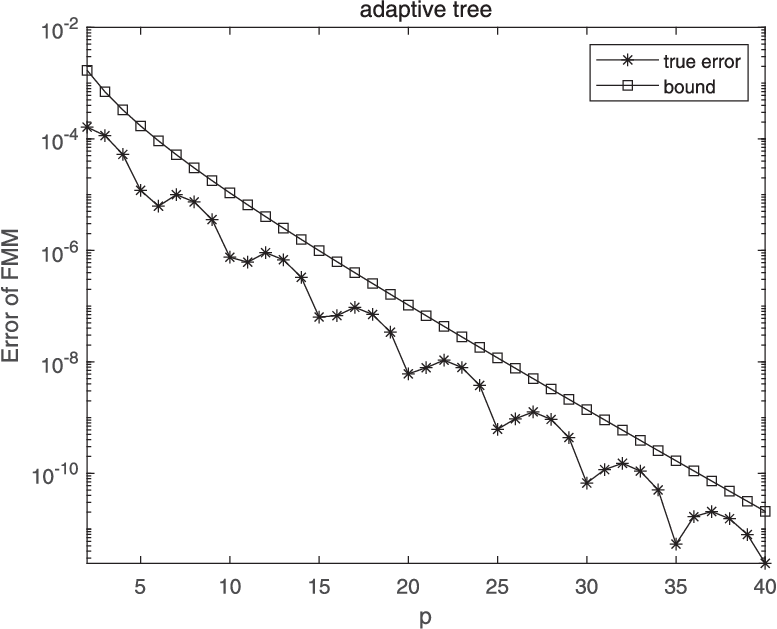}}\qquad\scalebox{0.55}{\includegraphics{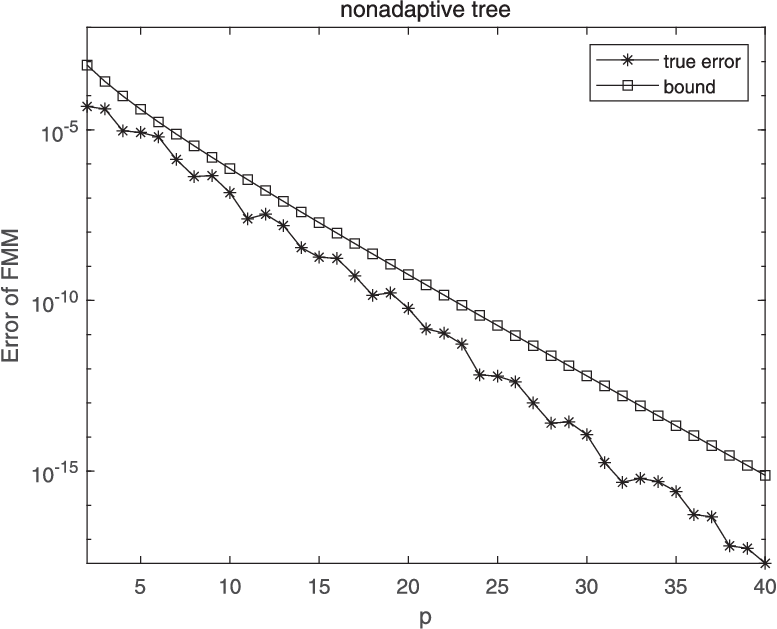}}\\
{\small Fig.4. The error of FMM for adaptive (left) and nonadaptive (right) trees.}
\end{figure}

In the above numerical examples, we only chose representative points in each tree structure. In fact, the error of FMM for the single-layer potential $(S\varphi)(\mathbf{x})$ is the sum of the errors at all field points.
Since only a few points (densely distributed points) in the adaptive tree structure meet $L_{\mathbf{y}_c}=L_{\mathbf{x}_c}-1$, while most points are similar to the nonadaptive one. Thus, when $p$ is not very large, the difference in the computational accuracy between the two tree structures is not as significant as expected.

\subsection{Error bound of the FMM for $(K\varphi)(\mathbf{x})$}
This section will study the bound on the error of the FMM for the double-layer potential $(K\varphi)(\mathbf{x})$. The discretized form of $(K\varphi)(\mathbf{x})$ is given by
$$(K\varphi)(\mathbf{x})=\frac{\mathrm{i}}{4}\sum_{j=1}^{N}\int_{\Delta\partial\Omega_j}\frac{\partial{H}^{(1)}_0(k|\mathbf{x}-\mathbf{y}_j|)}{\partial\nu(\mathbf{y}_j)}\varphi(\mathbf{y}_j)\mathrm{d}s(\mathbf{y}_j).$$
Since most of the expansions and translations about $(K\varphi)(\mathbf{x})$ are similar to those of $(S\varphi)(\mathbf{x})$, we only show the differences.

The multipole expansion:
$$\int_{\Delta\partial\Omega}\frac{\partial{H}^{(1)}_0(k|\mathbf{x}-\mathbf{y}|)}{\partial\nu(\mathbf{y})}\varphi(\mathbf{y})\mathrm{d}s(\mathbf{y})
=\sum_{n=-p}^p\mathcal{H}^+_n(k(\mathbf{x}-\mathbf{y}_c))\mathbf{M}_{n}(\mathbf{y}_c)+E_{ME},$$
where
\begin{equation}
\mathbf{M}_n(\mathbf{y}_c)=\int_{\Delta\partial\Omega}\frac{\partial\mathcal{J}^{-}_n(k(\mathbf{y}-\mathbf{y}_c))}{\partial\nu(\mathbf{y})}\varphi(\mathbf{y})\mathrm{d}s(\mathbf{y}),
\end{equation}
and
\begin{equation}
E_{ME}=\int_{\Delta\partial\Omega}\frac{\partial{R}_{0,p}^{H}\big(k(\mathbf{x}-\mathbf{y}_c),k(\mathbf{y}-\mathbf{y}_c)\big)}{\partial\nu(\mathbf{y})}\varphi(\mathbf{y})\mathrm{d}s(\mathbf{y}).
\end{equation}
The error of multipole moment $\mathbf{M}_n(\mathbf{y}_c)$:
\begin{equation}
EM_n=\int_{\Delta\partial\Omega}\frac{\partial{R}_{n,p}^{J}\big(k(\mathbf{y}_{c'}-\mathbf{y}_c),-k(\mathbf{y}-\mathbf{y}_{c'})\big)}{\partial\nu(\mathbf{y})}\varphi(\mathbf{y})\mathrm{d}s(\mathbf{y}).
\end{equation}
Except for these, the formulas of local expansion, M2M translation, M2L translation, L2L translation and $E_{MM}, E_{ML}, E_{LL}$ are the same as those of $(S\varphi)(\mathbf{x})$.

By the following recurrence relations of the Bessel functions
$$2\mathscr{B}'_n(z)=\mathscr{B}_{n-1}(z)-\mathscr{B}_{n+1}(z),\qquad \frac{2n}{z}\mathscr{B}_{n}(z)=\mathscr{B}_{n-1}(z)+\mathscr{B}_{n+1}(z),$$
we derive
$$\frac{\partial\mathscr{B}_n(k|\mathbf{y}|)e^{\pm{\mathrm{i}n}\theta_{\mathbf{y}}}}{\partial\nu(\mathbf{y})}=\frac{k}{2}\Big[\mathscr{B}_{n-1}
(k|\mathbf{y}|)e^{\pm\mathrm{i}\vartheta}-\mathscr{B}_{n+1}(k|\mathbf{y}|)e^{\mp\mathrm{i}\vartheta}\Big]e^{\pm{\mathrm{i}n}\theta_{\mathbf{y}}},$$
where $\vartheta$ is the angle between the vector $\mathbf{y}$ and the outward normal $\nu(\mathbf{y})$. It follows that
\begin{equation}
\left|\frac{\partial\mathscr{B}_n(k|\mathbf{y}|)e^{\pm{\mathrm{i}n}\theta_{\mathbf{y}}}}{\partial\nu(\mathbf{y})}\right|\leq\frac{k}{2}\Big[|\mathscr{B}_{n-1}(k|\mathbf{y}|)|+|\mathscr{B}_{n+1}(k|\mathbf{y}|)|\Big].
\end{equation}
Moreover, we have
\begin{eqnarray}
\left|\frac{\partial{R}_{m,p}^{\mathscr{B}}\big(k\mathbf{x},k\mathbf{y}\big)}{\partial\nu(\mathbf{y})}\right|\hspace{-0.6cm}&&\leq\left(\sum_{n=p+1}^\infty+\sum_{n=-\infty}^{-p-1}\right)
\left|\mathscr{B}_{m+n}(k|\mathbf{x}|)\right|\left|\frac{\partial{J}_n(k|\mathbf{y}|)e^{\mp{\mathrm{i}n}\theta_{\mathbf{y}}}}{\partial\nu(\mathbf{y})}\right|\nonumber\\
&&\leq\frac{k}{2}\sum_{n=p+1}^\infty\big[\left|\mathscr{B}_{n+m}(k|\mathbf{x}|)\right|+\left|\mathscr{B}_{n-m}(k|\mathbf{x}|)\right|\big]\big[\left|J_{n-1}(k|\mathbf{y}|)\right|+\left|J_{n+1}(k|\mathbf{y}|)\right|\big].
\end{eqnarray}
By (18), when $n\geq k|\mathbf{y}|$,
\begin{eqnarray}
\left|J_{n-1}(k|\mathbf{y}|)\right|+\left|J_{n+1}(k|\mathbf{y}|)\right|\hspace{-0.6cm}&&\leq\frac{1}{\Gamma(n)}\left(\frac{k|\mathbf{y}|}{2}\right)^{n-1}+\frac{1}{\Gamma(n+2)}\left(\frac{k|\mathbf{y}|}{2}\right)^{n+1}\nonumber\\
&&=\frac{1}{\Gamma(n)}\left(\frac{k|\mathbf{y}|}{2}\right)^{n-1}\left(1+\frac{k^2|\mathbf{y}|^2}{4(n+1)n}\right)\nonumber\\
&&\leq\frac{5}{4\Gamma(n)}\left(\frac{k|\mathbf{y}|}{2}\right)^{n-1}.
\end{eqnarray}
From (43), (44) and the proofs of Theorems 1 and 2, we derive
\begin{equation}
\left|\frac{\partial{R}_{m,p}^J\big(k\mathbf{x},k\mathbf{y}\big)}{\partial\nu(\mathbf{y})}\right|\leq\frac{5k}{2\Gamma(p+1)}\left(\frac{k|\mathbf{y}|}{2}\right)^p,
\end{equation}
and
\begin{equation}
\left|\frac{\partial{R}_{m,p}^H\big(k\mathbf{x},k\mathbf{y}\big)}{\partial\nu(\mathbf{y})}\right|\leq\frac{5kC_{p+m+1}(k|\mathbf{x}|)(2p+m-1)^{m+1}r^{p-1}}{\sqrt{2}\pi(k|\mathbf{x}|)^m(-(2p+m-1)\ln{r}-2m)}.
\end{equation}

Now, we can give the bounds on truncation errors of the FMM for double-layer potential $(K\varphi)(\mathbf{x})$ as follows.

{\theorem Given the conditions of Theorem $3$, the bounds on truncation errors of the FMM for $(K\varphi)(\mathbf{x})$ are given by}
\begin{eqnarray*}
&&|E_{ME}|\leq\frac{5kAC_{p+1}(k|\mathbf{x}-\mathbf{y}_c|)\mu^{p-1}}{-\sqrt{2}\pi\ln{\mu}}=\mathcal{O}\left(\mu^p\right),\\
&&|E_{MM}|\leq\frac{6\sqrt{2}kAC_p(k|\mathbf{x}-\mathbf{y}_c|)\eta^p}{\pi}=\mathcal{O}\left(\eta^p\right),\\
&&|E_{ML}|\leq\frac{15\sqrt{2}AC_{p+1}(k|\mathbf{x}_c-\mathbf{y}_c|)\lambda^p}{-\pi(3\ln\gamma+2)|\mathbf{x}_c-\mathbf{y}_c|}=\mathcal{O}\left(\lambda^p\right),\\
&&|E_{LL}|\leq\frac{16\sqrt{2}kAC_{p+1}(k\tau)\nu^{p+1}}{\pi}=\mathcal{O}\left(\nu^p\right).
\end{eqnarray*}

Theorem 4 can be derived by (42), (44), (45) and (46). The proof is quite similar to that of Theorem 3 and so is omitted.

\section{Conclusion and discussion}
This paper focuses on the estimation of the FMM error for scattering problems in 2-D. The error considered in this paper is the superposition of expansion errors and translation errors, in other words, the error of a translation chain: source$\rightarrow$multipole$\rightarrow$multipole$\rightarrow$local$\rightarrow$local. A novel bound on the FMM error and its convergence order are derived. We apply the result to the FMM with square quad-tree structure, and derive the specific error bound and convergence order.

Our error estimation has the following advantages: We consider the superposition of expansion errors and translation errors, which is closer to the real error of the algorithm. Our bound is sharper than the previous results and can be better used to the tree structure. The explicit formula of the error bound is given, and it does not contain any unknown constants.

The limiting forms of Bessel functions are used to study the error bound of the FMM in this paper, which have the similar forms with the functions used in the expansions and translations of the FMM for potential problems. Thus, the technique of this paper can be easily applied to the FMM for potential problems, and derive a more sharp bound on the FMM error. In addition, the proposed method and technique can also be applied to study the errors of FMM for other problems, such as elastostatic problems, Stokes flow problems.

\section*{Acknowledgements}
This work is supported by the National Natural Science Foundation of China (11201373) and Natural Science Foundation of Shaanxi Provincial Department of Education (14JK1747).

\end{document}